
\def\author{E. Esteves and S. Kleiman}
\def\title{The compactified Picard scheme\break
  of the compactified Jacobian}
\def\abstract{
 Let $C$ be an integral projective curve in any characteristic.  Given
an invertible sheaf $\cL$ on $C$ of degree 1, form the corresponding
Abel map $A_\cL\:C\to\CJ$, which maps $C$ into its compactified
Jacobian, and form its pullback map $\smash{A_\cL^*\:\Pic^0_\CJ\to J}$,
which carries the connected component of 0 in the Picard scheme back to
the Jacobian. If $C$ has, at worst, double points, then
$\smash{A_\cL^*}$ is known to be an isomorphism.  We prove that
$\smash{A_\cL^*}$ always extends to a map between the natural
compactifications, $\smash{\Pic^-_\CJ\to\CJ}$, and that the extended map
is an isomorphism if $C$ has, at worst, ordinary nodes and cusps.
 }

\def\PaperSize{letter} 
 \let\@=@  

\input amssym.def
\input amssym
\let\bt=\boxtimes

\def\GetNext#1 {\def\NextOne{#1}\if\relax\NextOne\let\next=\relax
        \else\let\next=\DoIt \fi \next}
\def\DoIt{\Act\NextOne\GetNext}
\def\ActOn#1{\expandafter\GetNext #1\relax\ }
\def\defcs#1{\expandafter\xdef\csname#1\endcsname}

\parskip=0pt plus 1.8pt \parindent10pt
\hsize 30pc 
\vsize 45pc 
\abovedisplayskip 4pt plus3pt minus1pt
\belowdisplayskip=\abovedisplayskip
\abovedisplayshortskip 2.5pt plus2pt minus1pt
\belowdisplayshortskip=\abovedisplayskip

\mag=\magstep1

\newskip\vadjustskip 
\def\centertext
 {\hoffset=\pgwidth \advance\hoffset-\hsize
  \advance\hoffset-2truein \divide\hoffset by 2\relax
  \voffset=\pgheight \advance\voffset-\vsize
  \advance\voffset-2truein \divide\voffset by 2\relax
  \advance\voffset\vadjustskip
 }
\newdimen\pgwidth\newdimen\pgheight
\def\letter{letter}\def\AFour{AFour}
\ifx\PaperSize\letter
 \pgwidth=8.5truein \pgheight=11truein
 \message{- Got a paper size of letter.  }\centertext
\fi
\ifx\PaperSize\AFour
 \pgwidth=210truemm \pgheight=297truemm
 \message{- Got a paper size of AFour.  }\centertext
\fi

\def\today{\ifcase\month\or     
 January\or February\or March\or April\or May\or June\or
 July\or August\or September\or October\or November\or December\fi
 \space\number\day, \number\year}
\footline={\hss\eightpoint\folio\hss}
 \newcount\pagenumber \pagenumber=1
 \def\advancepagenumber{\global\advance\pagenumber by 1}
\def\folio{\number\pageno} 
\headline={%
  \ifnum\pageno=0\hfill
  \else
   \ifnum\pageno=1\firstheadline
   \else
     \ifodd\pageno\oddheadline
     \else\evenheadline\fi
   \fi
  \fi
}
\expandafter\ifx\csname date\endcsname\relax \let\dato=\today
            \else\let\dato=\date\fi
\let\firstheadline\hfill
\def\oddheadline{\eightpoint
 \rlap{\dato} \hfil \headtitle
 }
\def\evenheadline{\eightpoint 
 \author\hfil\llap{\dato}}
\def\headtitle{\title}

 \font\twelvebf=cmbx12          
 \font\smc=cmcsc10              
 \font\tenbi=cmmi14
 \font\sevenbi=cmmi10 \font\fivebi=cmmi7
 \newfam\bifam  \textfont\bifam=\tenbi
 \scriptfont\bifam=\sevenbi \scriptscriptfont\bifam=\fivebi
 \mathchardef\variablemega="7121 
 \mathchardef\variablenu="7117 
\catcode`\@=11          
\def\eightpoint{\eightpointfonts
 \setbox\strutbox\hbox{\vrule height7\p@ depth2\p@ width\z@}%
 \eightpointparameters\eightpointfamilies
 \normalbaselines\rm
 }
\def\eightpointparameters{%
 \normalbaselineskip9\p@
 \abovedisplayskip9\p@ plus2.4\p@ minus6.2\p@
 \belowdisplayskip9\p@ plus2.4\p@ minus6.2\p@
 \abovedisplayshortskip\z@ plus2.4\p@
 \belowdisplayshortskip5.6\p@ plus2.4\p@ minus3.2\p@
 }
\newfam\smcfam
\def\eightpointfonts{%
 \font\eightrm=cmr8 \font\sixrm=cmr6
 \font\eightbf=cmbx8 \font\sixbf=cmbx6
 \font\eightit=cmti8
 \font\eightsmc=cmcsc8
 \font\eighti=cmmi8 \font\sixi=cmmi6
 \font\eightsy=cmsy8 \font\sixsy=cmsy6
 \font\eightsl=cmsl8 \font\eighttt=cmtt8}
\def\eightpointfamilies{%
 \textfont\z@\eightrm \scriptfont\z@\sixrm  \scriptscriptfont\z@\fiverm
 \textfont\@ne\eighti \scriptfont\@ne\sixi  \scriptscriptfont\@ne\fivei
 \textfont\tw@\eightsy \scriptfont\tw@\sixsy \scriptscriptfont\tw@\fivesy
 \textfont\thr@@\tenex \scriptfont\thr@@\tenex\scriptscriptfont\thr@@\tenex
 \textfont\itfam\eightit        \def\it{\fam\itfam\eightit}%
 \textfont\slfam\eightsl        \def\sl{\fam\slfam\eightsl}%
 \textfont\ttfam\eighttt        \def\tt{\fam\ttfam\eighttt}%
 \textfont\smcfam\eightsmc      \def\smc{\fam\smcfam\eightsmc}%
 \textfont\bffam\eightbf \scriptfont\bffam\sixbf
   \scriptscriptfont\bffam\fivebf       \def\bf{\fam\bffam\eightbf}%
 \def\rm{\fam0\eightrm}%
 }
\def\vfootnote#1{\insert\footins\bgroup
 \eightpoint\catcode`\^^M=5\leftskip=0pt\rightskip=\leftskip
 \interlinepenalty\interfootnotelinepenalty
  \splittopskip\ht\strutbox 
  \splitmaxdepth\dp\strutbox \floatingpenalty\@MM
  \leftskip\z@skip \rightskip\z@skip \spaceskip\z@skip \xspaceskip\z@skip
  \textindent{#1}\footstrut\futurelet\next\fo@t}

\def\p.{p.\penalty\@M \thinspace}
\def\pp.{pp.\penalty\@M \thinspace}
\newcount\sctno \newskip\sctnskip \sctnskip=0pt plus\baselineskip
\def\sctn#1\par
  {\removelastskip\vskip\sctnskip 
  \vskip-\sctnskip \bigskip\medskip
  \centerline{#1}\nobreak\medskip
}

\def\sct#1 {\sctno=#1\relax\sctn#1. }

\def\item#1 {\par\indent\indent\indent
 \hangindent3\parindent
 \llap{\rm (#1)\enspace}\ignorespaces}
 \def\inpart#1 {{\rm (#1)\enspace}\ignorespaces}
 \def\part {\par\inpart}

\def\Cs#1){\(\number\sctno.#1)}
\def\part#1 {\par\(#1)\enspace\ignorespaces}

\def\dsc#1 #2.{\medbreak\Cs#1) {\it #2.} \ignorespaces}
\def\proclaim#1 #2 {\medbreak
  {\bf#1 (\number\sctno.#2).  }\ignorespaces\bgroup\it}
\def\endproclaim{\par\egroup\medskip}
\def\pf{\endproclaim{\bf Proof.} \ignorespaces}
\def\lem{\proclaim Lemma } \def\prp{\proclaim Proposition }
  \def\thm{\proclaim Theorem }
\def\dfn#1 {\medbreak {\bf Definition (\number\sctno.#1).  }\ignorespaces}
\def\rmk#1 {\medbreak {\bf Remark (\number\sctno.#1).  }\ignorespaces}
 \newcount\refno \refno=0        \def\NoKey{*!*}
 \def\MakeKey{\advance\refno by 1 \expandafter\xdef
  \csname\TheKey\endcsname{{\number\refno}}\NextKey}
 \def\NextKey#1 {\def\TheKey{#1}\ifx\TheKey\NoKey\let\next\relax
  \else\let\next\MakeKey \fi \next}
 \def\RefKeys #1\endRefKeys{\expandafter\NextKey #1 *!* }
 \def\SetRef#1 #2,{\hang\llap
  {[\csname#1\endcsname]\enspace}{\smc #2},}
 \newbox\keybox \setbox\keybox=\hbox{[25]\enspace}
 \newdimen\keyindent \keyindent=\wd\keybox
\def\references{\kern-\medskipamount
  \sctn References\par
  \vskip-\medskipamount
  \bgroup   \frenchspacing   \eightpoint
   \parindent=\keyindent  \parskip=\smallskipamount
   \everypar={\SetRef}\par}
\def\endreferences{\egroup}

 \def\serial#1#2{\expandafter\def\csname#1\endcsname ##1 ##2 ##3
        {\unskip\ {\it #2\/} {\bf##1} (##2), ##3.}} 

\def\UThin{\penalty\@M \thinspace\ignorespaces}

\def\(#1){{\let~=\UThin\rm(#1)}}
\def\relaxnext@{\let\next\relax}
\def\cite#1{\relaxnext@
 \def\nextiii@##1,##2\end@{\unskip\space{\rm[\SetKey{##1},\let~=\UThin##2]}}%
 \in@,{#1}\ifin@\def\next{\nextiii@#1\end@}\else
 \def\next{{\rm[\SetKey{#1}]}}\fi\next}
\newif\ifin@
\def\in@#1#2{\def\in@@##1#1##2##3\in@@
 {\ifx\in@##2\in@false\else\in@true\fi}%
 \in@@#2#1\in@\in@@}
\def\SetKey#1{{\bf\csname#1\endcsname}}

\catcode`\@=12  

\let\:=\colon \let\ox=\otimes \let\x=\times
\let\wt=\widetilde \let\wh=\widehat

\let\into=\hookrightarrow \def\onto{\to\mathrel{\mkern-15mu}\to}
\let\To=\longrightarrow \def\TO#1{\buildrel#1\over\To}
\def\smashedlongrightarrow{\setbox0=\hbox{$\longrightarrow$}\ht0=1pt\box0}
\def\risom{\buildrel\sim\over{\smashedlongrightarrow}}
 \def\lgto{-\mathrel{\mkern-10mu}\to}
 \def\smashedlgto{\setbox0=\hbox{$\scriptstyle\lgto$}\ht0=1.85pt
        \lower1.25pt\box0}
\def\tto{\buildrel\lgto\over{\smashedlgto}}
\let\vf=\varphi    \let\?=\overline
\let\ve=\varepsilon

\def\CJ{{\bar J}}  

 \def\IP{{\bf P}}   
\def\IF{{\bf F}}
\def\Act#1{\defcs{c#1}{{\cal#1}}}
 \ActOn{A B C D E F G H I J K L M N O P Q R T U V X Y }
\def\Act#1{\defcs{#1}{\mathop{\rm#1}\nolimits}}
 \ActOn{Ass cod Div Ext Hilb length Pic Quot sm Spec Supp }
\def\Act#1{\defcs{c#1}{\mathop{\it#1}\nolimits}}
 \ActOn{Cok Ext Hom Ker Sym }
\def\Act#1{\defcs{#1}{\hbox{\rm\ #1 }}}
 \ActOn{and by for where with on }

\catcode`\@=11

 \def\activeat#1{\csname @#1\endcsname}
 \def\def@#1{\expandafter\def\csname @#1\endcsname}
 {\catcode`\@=\active \gdef@{\activeat}}

\let\ssize\scriptstyle
\newdimen\ex@   \ex@.2326ex

 \def\requalfill{\cleaders\hbox{$\mkern-2mu\mathord=\mkern-2mu$}\hfill
  \mkern-6mu\mathord=$}
 \def\eqfill{$\m@th\mathord=\mkern-6mu\requalfill}
 \def\deffill{\hbox{$:=$}$\m@th\mkern-6mu\requalfill}
 \def\fiberbox{\hbox{$\vcenter{\hrule\hbox{\vrule\kern1ex
     \vbox{\kern1.2ex}\vrule}\hrule}$}}
 \def\Fiberbox{\rlap{\kern-0.75pt\raise0.75pt\fiberbox}%
        \kern.75pt{\lower0.75pt\fiberbox}}

 \font\arrfont=line10
 \def\Swarrow{\vcenter{\hbox{$\swarrow$\kern-.26ex
    \raise1.5ex\hbox{\arrfont\char'000}}}}

 \newdimen\arrwd
 \newdimen\minCDarrwd 
   \setbox\z@\hbox{$\rightarrow$} \minCDarrwd=\wd\z@
        
 \def\findarrwd#1#2#3{\arrwd=#3%
  \setbox\z@\hbox{$\ssize\;{#1}\;\;$}%
 \setbox\@ne\hbox{$\ssize\;{#2}\;\;$}%
  \ifdim\wd\z@>\arrwd \arrwd=\wd\z@\fi
  \ifdim\wd\@ne>\arrwd \arrwd=\wd\@ne\fi}
 \newdimen\arrowsp\arrowsp=0.375em
 \def\findCDarrwd#1#2{\findarrwd{#1}{#2}{\minCDarrwd}
    \advance\arrwd by 2\arrowsp}
 \newdimen\minarrwd 
 \setbox\z@\hbox{$\rightarrow$} \minarrwd=\wd\z@

 \def\harrow#1#2#3#4{{\minarrwd=#1\minarrwd%
   \findarrwd{#2}{#3}{\minarrwd}\kern\arrowsp
    \mathrel{\mathop{\hbox to\arrwd{#4}}\limits^{#2}_{#3}}\kern\arrowsp}}
 \def@]#1>#2>#3>{\harrow{#1}{#2}{#3}\rightarrowfill}
 \def@>#1>#2>{\harrow1{#1}{#2}\rightarrowfill}
 \def@<#1<#2<{\harrow1{#1}{#2}\leftarrowfill}
 \def@={\harrow1{}{}\eqfill}
 \def@:#1={\harrow1{}{}\deffill}
 \def@ N#1N#2N{\vCDarrow{#1}{#2}\UpDownarrow}
 \def\UpDownarrow{\uparrow\,\Big\downarrow}

\def@'#1'#2'{\harrow1{#1}{#2}\tarrowfill}
 \def\lgTo{\dimen0=\arrwd \advance\dimen0-2\arrowsp
        \hbox to\dimen0{\rightarrowfill}}
 \def\smashedlgTo{\setbox0=\hbox{$\scriptstyle\lgTo$}\ht0=1.85pt
        \lower1.25pt\box0}
 \def\tto{\buildrel\lgTo\over{\smashedlgTo}}
 \def\tarrowfill{\hfil$\tto$\hfil}  

 \def@={\ifodd\row\harrow1{}{}\eqfill
   \else\vCDarrow{}{}\Vert\fi}
 \def@.{\ifodd\row\relax\harrow1{}{}\hfill
   \else\vCDarrow{}{}.\fi}
 \def@|{\vCDarrow{}{}\Vert}
 \def@ V#1V#2V{\vCDarrow{#1}{#2}\downarrow}
\def@ A#1A#2A{\vCDarrow{#1}{#2}\uparrow}
 \def@(#1){\arrwd=\csname col\the\col\endcsname\relax
   \hbox to 0pt{\hbox to \arrwd{\hss$\vcenter{\hbox{$#1$}}$\hss}\hss}}

 \def\squash#1{\setbox\z@=\hbox{$#1$}\finsm@@sh}
\def\finsm@@sh{\ifnum\row>1\ht\z@\z@\fi \dp\z@\z@ \box\z@}

 \newcount\row \newcount\col \newcount\numcol \newcount\arrspan
 \newdimen\vrtxhalfwd  \newbox\tempbox

 \def\innernewdimen{\alloc@1\dimen\dimendef\insc@unt}
 \def\measureinit{\col=1\vrtxhalfwd=0pt\arrspan=1\arrwd=0pt
   \setbox\tempbox=\hbox\bgroup$}
 \def\setinit{\col=1\hbox\bgroup$\ifodd\row
   \kern\csname col1\endcsname
   \kern-\csname row\the\row col1\endcsname\fi}
 \def\findvrtxhalfsum{$\egroup
  \expandafter\innernewdimen\csname row\the\row col\the\col\endcsname
  \global\csname row\the\row col\the\col\endcsname=\vrtxhalfwd
  \vrtxhalfwd=0.5\wd\tempbox
  \global\advance\csname row\the\row col\the\col\endcsname by \vrtxhalfwd
  \advance\arrwd by \csname row\the\row col\the\col\endcsname
  \divide\arrwd by \arrspan
  \loop\ifnum\col>\numcol \numcol=\col%
 \expandafter\innernewdimen \csname col\the\col\endcsname
     \global\csname col\the\col\endcsname=\arrwd
   \else \ifdim\arrwd >\csname col\the\col\endcsname
      \global\csname col\the\col\endcsname=\arrwd\fi\fi
   \advance\arrspan by -1 %
   \ifnum\arrspan>0 \repeat}
 \def\setCDarrow#1#2#3#4{\advance\col by 1 \arrspan=#1
    \arrwd= -\csname row\the\row col\the\col\endcsname\relax
    \loop\advance\arrwd by \csname col\the\col\endcsname
     \ifnum\arrspan>1 \advance\col by 1 \advance\arrspan by -1%
     \repeat
    \squash{\mathop{
     \hbox to\arrwd{\kern\arrowsp#4\kern\arrowsp}}\limits^{#2}_{#3}}}
 \def\measureCDarrow#1#2#3#4{\findvrtxhalfsum\advance\col by 1%
   \arrspan=#1\findCDarrwd{#2}{#3}%
    \setbox\tempbox=\hbox\bgroup$}


\def\vCDarrow#1#2#3{\kern\csname col\the\col\endcsname
    \hbox to 0pt{\hss$\vcenter{\llap{$\ssize#1$}}%
     \Big#3\vcenter{\rlap{$\ssize#2$}}$\hss}\advance\col by 1}

 \def\setCD{\def\harrow{\setCDarrow}%
  \def\\{$\egroup\advance\row by 1\setinit}
  \m@th\lineskip3\ex@\lineskiplimit3\ex@ \row=1\setinit}
 \def\endsetCD{$\strut\egroup}
 \def\measure{\bgroup
  \def\harrow{\measureCDarrow}%
  \def\\##1\\{\findvrtxhalfsum\advance\row by 2 \measureinit}%
  \row=1\numcol=0\measureinit}
 \def\endmeasure{\findvrtxhalfsum\egroup}

\newbox\CDbox \newdimen\sdim

 \newcount\savedcount
 \def\CD#1\endCD{\savedcount=\count11%
   \measure#1\endmeasure
   \vcenter{\setCD#1\endsetCD}%
   \global\count11=\savedcount}

 \catcode`\@=\active
 \RefKeys
 AIK76 AK76 AK79 AK80  E95 E97 EGK00 EGK02 FGS99 FGA EGA 
 K87 L02 
S04
 \endRefKeys
{\leftskip=0pt plus1fill \rightskip=\leftskip
 \obeylines
 \leavevmode \bigskip
 {\twelvebf \title
 } \medskip
 \footnote{}{\noindent %
 MSC-class: 14H40 (Primary) 14K30, 14H20  (Secondary).}
 Eduardo Esteves\footnote{$^{1}$}{%
    Supported in part by Programa do Mil\^enio -- AGIMB and CNPq
    Proc. 300004/95-8 (NV).}
 {\eightpoint\it\medskip
 Instituto de Matem\'atica Pura e Aplicada
 Estrada D. Castorina {\sl110, 22460--320} Rio de Janeiro RJ, BRAZIL
 \rm E-mail: \tt Esteves\@impa.br \medskip
 } and \medskip
 Steven Kleiman\footnote{$^{2}$}{%
    Supported in part by NSF grant 9400918-DMS
 and Programa do Mil\^enio -- AGIMB.}
 {\eightpoint\it\medskip
 Department of Mathematics, Room {\sl 2-278} MIT,
 {\sl77} Mass Ave, Cambridge, MA {\sl02139-4307}, USA
 \rm E-mail: \tt Kleiman\@math.mit.edu \medskip
 \rm \dato \bigskip
 }
}
{
 \advance\leftskip by 1.5\parindent \advance\rightskip by 1.5\parindent
 \eightpoint \noindent
 {\smc Abstract.}\enspace \ignorespaces \abstract \par}
\def\headtitle{The compactified Picard scheme
  of the compactified Jacobian}

\sct1 Introduction

Let $C$ be an integral projective curve of arithmetic genus $g$, defined
over an algebraically closed field of any characteristic.  Form its
(generalized) Jacobian $J$, the connected component of the identity of
the Picard scheme of $C$.  If $C$ is singular, then $J$ is not
projective.  So for about forty years, numerous authors have studied a
natural compactification of $J$: the (fine) moduli space $\CJ$ of
torsion-free sheaves of rank 1 and degree 0 on $C$.  It is called the
{\it compactified Jacobian}.

Recently, the compactified Jacobian appeared in Laumon's preprint
\cite{L02}, where he identified, up to homeomorphism, affine Springer
fibers with coverings of compactified Jacobians.  For that
identification, he used the autoduality of the compactified Jacobian, a
property established in \cite{EGK02} and explained next.

{}From now on, assume $C$ has, at worst, points of multiplicity 2 (or
double points).  
For each invertible sheaf $\cL$ of degree 1 on $C$, form the Abel map
$A_\cL\:C\to\CJ$, given by $P\mapsto\cM_P\ox\cL$ where $\cM_P$ is the
ideal sheaf of $P$; it is a closed embedding if $C$ is not of genus 0.
Form the pullback map
        $$A_\cL^*\:\Pic^0_\CJ\to J,$$
 carrying the connected component of 0 in the Picard scheme back to the
Jacobian.  Then $A_\cL^*$ is an isomorphism and is independent of $\cL$;
see \cite{EGK02, Thm.~2.1, \p.595}.

Since the singularities are locally planar, $\CJ$ is integral by
\cite{AIK76, (9), \p.8}.  Hence, not only does $\Pic^0_\CJ$ exist, but
also it admits a natural compactification: its closure $\Pic^-_\CJ$ in
the compactified Picard scheme $\Pic^=_\CJ$, the (fine) moduli space of
torsion-free sheaves of rank 1 on $\CJ$; see \cite{AK79, Thm.~(3.1),
\p.28}.  Does $A_\cL^*$ extend to a map between the compactifications?
If so, then is the extension an isomorphism?

These  questions were posed to the authors by Sawon.  As
mentioned in his introduction to \cite{S04}, his results on dual 
fibrations to fibrations by Abelian varieties, in the ``nicest'' 
cases, depend on ``extending autoduality to the 
compactifications.''

It is not true, for every map, that the pullback of a torsion-free sheaf
is still torsion free.  But, for $A_\cL$, it is true!  There are two
basic reasons why: first, $A_\cL^*$ is independent of $\cL$; second, the
maps $A_\cL$ can be bundled up into a {\it smooth} map $C\x J\to\CJ$.
Thus there exists an extended pullback map
	$$A_\cL^*\:\Pic^-_\CJ\to \CJ;$$
this existence statement is the content of Theorem (2.6) below.  (The
statement and its proof already appear on the web in the preliminary
version of \cite{EGK02}.)

Is the extended map $A_\cL^*$ also an isomorphism?  This question seems
much harder.

{}From now on, assume that $C$ has, at worst, ordinary nodes and cusps.
Then the extended map $A_\cL^*$ is, indeed, an isomorphism, according to
Theorem (4.1), our main result.  Here is a sketch of the proof.

First, recall the definition of the inverse $\beta\:J\to\Pic^0_\CJ$ from
\cite{EGK02, Prop.~2.2, \p.595}, or rather \cite{EGK02, Rmk.~2.4,
\p.597}.  Let $\cI$ be the universal sheaf on $C\x\CJ$, and $\cP$ the
determinant of cohomology of $\cI\ox\cL^{\ox g-1}$ with respect to the
projection $C\x\CJ\to\CJ$.  The sheaf $\cP^{-1}$ has a canonical regular
section, whose zero scheme we denote by $\Theta$ and call the {\it theta
divisor\/} associated to $\cL$.  Let $p_1\:\CJ\x J\to\CJ$ be the
projection, and $\mu\:\CJ\x J\to\CJ$ the multiplication map.  Form the
invertible sheaf
	$$\cT:=\cO_{\CJ\x J}(\mu^*\Theta-p_1^*\Theta).$$ 
 Then $\cT$ defines the map $\beta\: J\to\Pic^0_\CJ$.

We need only show that $\beta$ extends to a map
$\?\beta\:\CJ\to\Pic^-_\CJ$, or that $\cT$ extends to a
sheaf on $\CJ\x\CJ$ that is flat for $p_2\:\CJ\x\CJ\to\CJ$ and whose
fibers are torsion-free  rank-1.  In fact, as $\?\beta$ is
unique if it exists, we need only find a flat surjection
$\zeta\:W\to\CJ$ such that the pullback $\cT_\zeta$ of $\cT$ to
$\CJ\x(\zeta^{-1}J)$ extends to a torsion-free  $\cH$ on
$\CJ\x W/W$.

To define $\zeta$, let $C_0$ be the smooth locus, and $P_1,\dots,P_m$
the singular points of $C$.  Set $C_i:=C_0\cup\{P_i\}$ for $1\le i\le
m$.  Set
	$$H_m:=C_m\x\cdots\x C_1\subseteq C^{\x m}
	\and W:=H_m\x J.$$
Define $\zeta\:W\to\CJ$ as the natural map given
on a pair, consisting of an $m$-tuple of points $(Q_m,\dots,Q_1)$ of $C$
and an invertible sheaf $\cN$ on $C$ of degree 0, by
	$$\zeta(Q_m,\dots,Q_1,\cN):=\cM_{Q_m}\ox\cdots\ox\cM_{Q_1}\ox
	\cN\ox\cL^{\ox m}$$
where $\cM_{Q_i}$ is the ideal of $Q_i$ in $C$ for $i=1,\dots, m$.
Since $\zeta$ is a composition of base extensions of bigraded Abel maps,
$\zeta$ is smooth by \cite{EGK00, Cor.~2.6, \p.5969}.  It is also
surjective, owing to our assumption on the singularities: since, at each
local ring $\cO_{C,P_i}$, the singularity degree $\delta$ is 1, any
torsion-free rank-1 module is either free or isomorphic to the maximal
ideal.

Next, we resolve the rational map $\mu\circ(1\x\zeta)\:\CJ\x W\to\CJ$
by using the $n$-{\it flag schemes} $F_n$ of $\cI/(C\x\CJ)/\CJ$.  The scheme
$F_n$ is the (fine) moduli space parameterizing $n$-chains of
torsion rank-1 sheaves on $C$:
	$$\cI_n\subset\cI_{n-1}\subset\cdots\subset\cI_1
	\subset\cI_0$$
 where $\cI_0$ is of degree $0$ and where each quotient $\cI_{i-1}/\cI_i$ 
is of length 1.

The scheme $F_n$ comes equipped with two important maps.  The first is
the {\it multiplication map} $\gamma_n\:F_n\x J\to\CJ$, which sends a pair
consisting of a chain as above and an invertible sheaf $\cN$ on $C$ of
degree 0 to the tensor product $\cI_n\ox\cN\ox\cL^{\ox n}$.   The second
is the {\it resolution map}  $\wh\psi^{(n)}\:F_n\to \CJ\x C^{\x n}$, 
 which sends a chain as above to the pair consisting of the 
sheaf $\cI_0$ and the $n$-tuple $(Q_n,\dots,Q_1)\in C^{\x n}$ such 
that $\cI_{i-1}/\cI_i$ is supported on $Q_i$ for each $i$.

Set $\wt F_m:=(\wh\psi^{(m)})^{-1}(\CJ\x H_m)$.  Then we have the
following diagram, in which the right vertical map $\mu$ is only
rational:
	$$\CD
           \wt F_m\x J        @>\gamma_m>>       \CJ        \\
        @V\wh\psi^{(m)}\x 1VV                       @A\mu AA        \\
        \CJ\x H_m\x J      @>1\x\zeta>>    \CJ\x\CJ
	\endCD$$
 Note that $\zeta^{-1}(J)=C_0^{\x m}\x J$.  It is not difficult to see
that $\wh\psi^{(m)}$ restricts to an isomorphism over $\CJ\x C_0^{\x m}$.
Also, the composition $\mu\circ(1\x\zeta)\circ(\wh\psi^{(m)}\x 1)$ is
defined on $(\wh\psi^{(m)}\x 1)^{-1}(\CJ\x C_0^{\x m})$ and agrees with
$\gamma_m$.  Therefore, $\cT_\zeta$ extends to the sheaf
	$$\cH:=\cO_{\CJ\x W}(-q_1^*\Theta)\ox
	(\wh\psi^{(m)}\x 1)_*\cO_{\wt F_m\x J}(\gamma_m^*\Theta)$$
 where $q_1\:\CJ\x W\to\CJ$ is the projection.

The delicate part is now to prove that $\cH$ is flat over $W$ with
torsion-free rank-1 fibers.  It suffices to prove that $(\zeta\x1)^*\cH$
on $W\x W$ is flat over $W$ with torsion-free rank-1 fibers.  Now, there
are base-change formulas for the determinant of cohomology and for the
direct image.  As the corresponding base-change maps, we use the
horizontal maps in the following natural Cartesian square:
	$$\CD
 \wt F_{2m}\x J		@>>>	\wt F_m\x J		\\
     @V\xi VV				@VV\wh\psi^{(m)}\x1V	\\
        W\x W		@>\zeta\x1>>	\CJ\x W
	\endCD$$
where $\wt F_{2m}:=(\wh\psi^{(2m)})^{-1}(J\x H_m\x H_m)$ and $\xi$ is,
up to switching factors, $\wh\psi^{(2m)}\x 1$.  We now apply the
technical, but essential, Lemma (3.3) to complete the proof.

All our results are, in fact, proved not just for an individual curve
$C$, but for a flat projective family of
 (geometrically)
 integral curves over an arbitrary base scheme.  All schemes are
implicitly assumed to be locally Noetherian.

In short, in Section 2, we prove that the autoduality map $A_\cL^*$
extends if the curves have, at worst, double points.  In Section 3, we
study flag schemes.  Finally, in Section 4, we prove that the extended
map $A_\cL^*$ is an isomorphism if the curves have, at worst, ordinary
nodes and cusps.

\sct2 Extension

\dsc1 The compactified Jacobian.  By a {\it flat projective family of
integral curves} $C/S$, let us mean that $C$ is a flat and projective
$S$-scheme with geometrically integral fibers of dimension 1.

Given such a family $C/S$ and given an integer $n$, recall from
\cite{AK76}, \cite{AK79}, \cite{AK80}, and \cite{EGK02}
that there exists a projective $S$-scheme $\CJ^n$, or $\CJ^n_{C/S}$,
parameterizing the torsion-free rank-1 sheaves of degree $n$ on the
fibers of $C/S$.  And there exists an open subscheme $J^n$, or
$J^n_{C/S}$, parameterizing those sheaves that are invertible.
Furthermore, forming these schemes commutes with changing the base $S$.
For short, set $J:=J^0$ and $\CJ:=\CJ^0$.  Customarily, $J$ is called
the (relative generalized) {\it Jacobian} of $C/S$, and $\CJ$ the {\it
compactified Jacobian}.

More precisely, $\CJ^n$ represents the \'etale sheaf associated to the
functor whose $T$-points are degree-$n$ torsion-free  rank-1 sheaves
$\cI$ on $C\x T/T$.  Such an $\cI$ is a $T$-flat coherent sheaf on $C\x
T$ such that, for each point $t$ of $T$, the fiber $\cI(t)$ is
torsion-free and of generic rank 1 on the fiber $C(t)$ and also
	$$\chi(\cI(t))-\chi(\cO_{C(t)})=n.$$

\dsc2 Multiplication and translation. 
Let $C/S$ be a flat projective family of integral curves.  Let $m$ and 
$n$ be arbitrary integers.  Let 
$U^{n,m}\subseteq\CJ^n\x\CJ^m$ be the open subscheme that 
represents the \'etale sheaf associated to the subfunctor whose 
$T$-points are the pairs of torsion-free rank-1 sheaves $(\cI,\cJ)$ 
on $C\x T/T$ such that $\cI$ is invertible where $\cJ$ is not.  
Define the  {\it multiplication} map
	$$\mu\: U^{n,m}\to\CJ^{m+n} \hbox{ by }
	\mu(\cI,\cJ):=\cI\ox\cJ.$$

For each invertible sheaf $\cM$ of degree $m$ on $C/S$, define the {\it
translation\/} by $\cM$
        $$\mu_\cM\:\CJ^n\to\CJ^{m+n} \hbox{ by }
	\mu_\cM(\cI):=\cI\ox\cM.$$

 In other words,  $\cM$ defines a section $\sigma\:S\to J^m$, and
$\mu_\cM:=\mu\circ(1\x\sigma)$; the composition makes sense
because $U^{n,m}\supseteq\CJ^n\x J^m$.

\dsc3 The (bigraded) Abel map.  Let $C/S$ be a flat projective family of
integral curves.  Let $\Delta\subset C\x C$ be the diagonal.  Then its
ideal defines a map $\iota\:C\to\CJ^{-1}$.

Let $m$ be an arbitrary integer.  Let $W^m\subseteq C\x\CJ^{m+1}$ be the
inverse image of $U^{-1,m+1}$ under $\iota\x 1$.  Define the {\it
bigraded Abel map} to be the composition
	$$A:=\mu\circ(\iota\x 1)\:W^m\to\CJ^m.$$

Let $\cL$ be an invertible sheaf of degree 1 on $C/S$, and $\sigma\:S\to
J^1$ be the corresponding section.  Define the {\it Abel map} associated
to $\cL$ to be the composition
        $$A_\cL:=A\circ(1\x\sigma)\:C\to \CJ;$$
the composition makes sense because $C\x J^1\subseteq W^0$.

\dsc4 Autoduality.  Let $C/S$ be a flat projective family of integral
curves.  Assume that the curves (the geometric fibers of $C/S$) are
locally planar.  Then the projective $S$-scheme $\CJ^n$ is flat, and its
geometric fibers are integral local complete intersections; see
\cite{AIK76, (9), \p.8}.  Hence, the Picard scheme $\Pic_{\CJ^n/S}$
exists, and is a disjoint union of quasi-projective $S$-schemes; see
Th\'eor\`eme~3.1, \p.232-06, in \cite{FGA}, and Corollary~(6.7)(ii),
\p.96, in \cite{AK80}.  Also, by \cite{AK79, Thm.~(3.1), \p.28}, there
exists an $S$-scheme $\smash{\Pic^=_{\CJ^n/S}}$, the {\it compactified Picard
scheme}, that parameterizes torsion-free rank-1 sheaves on the fibers of
$\CJ^n/S$; moreover, the connected components of $\Pic^=_{\CJ^n/S}$ are
proper over $S$.

As is customary \cite{FGA, \p.236-03}, let $\Pic^0_{\CJ/S}$ denote the
set-theoretic union of the connected components of the identity 0 in the
fibers of $\smash{\Pic_{\CJ/S}}$, and let $\smash{\Pic^\tau_{\CJ/S}}$
denote the set of points of $\Pic_{\CJ/S}$ that have a multiple in
$\Pic^0_{\CJ/S}$.  The set $\Pic^\tau_{\CJ/S}$ is open; give it the
induced scheme structure.  Then, by general principles, forming
$\Pic_{\CJ/S}$ and $\Pic^\tau_{\CJ/S}$ commutes with changing $S$.

Denote by $\smash{\Pic^-_{\CJ/S}}$ the schematic closure of
$\smash{\Pic^\tau_{\CJ/S}}$ in $\smash{\Pic^=_{\CJ^n/S}}$.  Note that
forming $\smash{\Pic^-_{\CJ/S}}$ commutes with changing $S$ via a flat
map: it does so topologically because a flat map is open; whence, it
does so schematically because a flat map carries associated points to
associated points.  If the fibers of $C/S$ have, at worst, nodes and
cups, then forming $\smash{\Pic^-_{\CJ/S}}$ commutes with changing $S$
via an arbitrary map, owing to Theorem (4.1), our main result, since
forming $\CJ$ does so.

The Abel map $A_\cL$ induces an $S$-map,
        $$\textstyle A_\cL^*\:\Pic_{\CJ/S}\to \coprod_nJ^n.$$
 By \cite{EGK02, Thm.~2.1, \p.595}, if the geometric fibers of 
$C/S$ have, at worst, double points, then 
$\Pic^0_{\CJ/S}=\Pic^\tau_{\CJ/S}$.  Furthermore, the Abel 
map induces an isomorphism,
        $$A_\cL^*\:\Pic^0_{\CJ/S}\risom J,$$
 which is independent of the choice of the invertible sheaf $\cL$ of
degree 1 on $C/S$; in fact, the isomorphism exists whether or not any
sheaf $\cL$ does.  Let us call this isomorphism the {\it autoduality
isomoprhism}.

 \prp5 Let $C/S$ be a flat projective family of integral curves, $m$ and
$n$ integers, $\cM$ an invertible sheaf of degree $m$ on $C/S$.  Suppose
the curves have, at worst, double points.  Then the translation map
$\mu_\cM$
induces an isomorphism
        $$\mu_\cM^*\:\Pic^0_{\CJ^{m+n}/S}\risom\Pic^0_{\CJ^n/S},$$
 which is independent of $\cM$.  If $m=0$, then $\mu_\cM^*$ is equal to
the identity.
 \pf Note that $\mu_{\cO_C}=1_{\CJ^n}$.  And, if $\cM_1$ is also an
invertible sheaf on $C$, then
        $$\mu_\cM\circ\mu_{\cM_1}=\mu_{\cM\ox\cM_1}.$$
  So $\mu_\cM$ is an isomorphism, whose inverse is 
$\mu_{\cM^{-1}}$. Hence $\mu_\cM^*$ is an isomorphism.  Moreover, 
if $\cM_1$ is of degree $m$ too, then $\cM\ox\cM_1^{-1}$ is of 
degree 0, and it suffices to prove that 
$\mu^*_{\cM\ox\cM_1^{-1}}=1.$   Thus we may assume $m=0$.

To prove that $\mu_\cM^*=1$, we may change the base via an \'etale
covering, and so assume that the smooth locus of $C/S$ admits a 
section $\sigma$.  Set $\cL:=\cO_C(\sigma(S))$.  Then $\cL$  is an 
invertible sheaf on $C$.  So,
        $$\mu_\cM=\mu_{\cL^{\ox n}}\circ\mu_\cM\circ
        \mu_{\cL^{\ox -n}}.$$
 Hence, since $\cL$ is of degree 1 on $C/S$, we may assume $n=0$.

Note that $\mu_\cM\circ A_\cL=A_{\cM\ox\cL}$.  Now,
$A_{\cM\ox\cL}^*=A_\cL^*$ and $A_\cL^*$ is an isomorphism; see (2.4).
Thus $\mu_\cM^*=1$, and the proof is complete.

 \thm6 Let $C/S$ be a flat projective family of integral curves with, at
worst, double points.  Then the autoduality isomorphism
$\smash{\Pic^0_{\CJ/S}}\risom J$ extends uniquely to a map of
compactifications $\smash{\Pic^-_{\CJ/S}}\to\CJ$.
 \pf Set $U:=\Pic^0_{\CJ/S}$ and $\?U:=\Pic^-_{\CJ/S}$.  Since $J/S$ is
smooth and admits a section (for example, the 0-section), by \cite{AK79,
Thm.~(3.4)(iii), \p.40}, $\CJ\x\?U\bigm/\?U$ carries a universal sheaf
$\cP$, which is determined up to tensor product with the pullback of an
invertible sheaf on $\?U$.

The extension $\eta\:\?U\to\CJ$ of the autoduality isomorphism is unique
if it exists, because $U$ is schematically dense in $\?U$ and $\CJ$ is
separated.  Hence, by descent theory, it suffices to construct $\eta$
after changing the base via an \'etale covering.  So we may assume that
the smooth locus of $C/S$ admits a section $\sigma$.  Set
$\cL:=\cO_C(\sigma(S))$.  Then $\cL$ is invertible of degree 1 on $C/S$.
So the autoduality isomorphism is simply $A_\cL^*$, and it suffices to
prove that $(A_\cL\x 1)^*\cP$ is torsion-free rank-1 on
$C\x\?U\bigm/\?U$.

Form the bigraded Abel map $A\:C\x J^1\to\CJ$.  It is smooth by
\cite{EGK00, Cor.~2.6, \p.5969}.  Hence $(A\x 1)^*\cP$ is
torsion-free rank-1 on $C\x J^1\x\?U\bigm/\?U$.  It suffices to prove
that $(A\x1)^*\cP$ is torsion-free rank-1 on 
$C\x J^1\x\?U\big/(J^1\x\?U)$, since $(A_\cL\x1)^*\cP$ is its fiber
over the point of $J^1$ representing $\cL$.  Now, $(A\x1)^*\cP$ is
flat over $J^1\x\?U$, by the local criterion, if its fiber is flat over
the fiber $J^1(u)$ for each $u\in\?U$.

  Fix a $u\in\?U$.  Making a suitable faithfully flat base change
$S'/S$, we may assume that the field $k(u)$ is equal to the field of the
image of $u$ in $S$.  Set $\cI:=\cP(u)$.  It suffices to prove that
$A(u)^*\cI$ is a torsion-free rank-1 sheaf on $C(u)\x
J^1(u)\big/J^1(u)$.

Let  $\cM$ be an invertible sheaf of degree 0 on $C/S$.  Then 
the translation map $\mu_\cM$ gives rise to the following 
commutative diagram:
        $$\CD
        C\x J^1 @>A>> \CJ \\
        @V1\x\mu_\cM VV        @V\mu_\cM VV \\
        C\x J^1 @>A>> \CJ
        \endCD$$
 By Proposition~\Cs5), $\mu_\cM^*$ is the identity on 
$\Pic^0_{\CJ/S}$, so on its closure $\?U$ too.  Thus
$\mu_\cM(u)^*\cI=\cI$.  Now, the diagram is commutative; hence,
        $$(1\x\mu_\cM(u))^*A(u)^*\cI=A(u)^*\cI.\eqno\Cs6.1)$$

Since $J^1(u)$ is integral, the lemma of generic flatness applies, and
it implies that there is a dense open subset $W$ of $J^1(u)$ over which
$A(u)^*\cI$ is flat.  Now, by Part~(ii)(a) of Lemma (5.12) on \p.85 of
\cite{AK80}, it is an open condition on the base for a flat family of
sheaves to be torsion-free rank-1 provided they are supported on a
family whose geometric fibers are integral of the same dimension.
Hence, since $A(u)^*\cI$ is torsion-free rank-1 on $C(u)\x
J^1(u)\big/k(u)$, after shrinking $W$, we may assume that the
restriction of $A(u)^*\cI$ to $C\x W\bigm/W$ is torsion-free rank-1.  Fix
an arbitrary point $j_1$ of $W$ and one $j_2$ of $J^1(u)$.

Making a suitable faithfully flat base change $S'/S$, we may assume 
that each of $j_1$ and $j_2$ lies in the image of a section of 
$J^1/S$.  These sections represent invertible sheaves $\cM_1$ and 
$\cM_2$ of degree 1 on $C/S$; set $\cM:=\cM_1\ox \cM_2^{-1}$.  
Equation \Cs6.1) implies that $A(u)^*\cI$ is torsion-free rank-1 
over $\mu_\cM(u)^{-1}W$ as well.  Now, $j_2$ is an arbitrary point 
of $J^1(u)$.  Hence $A(u)^*\cI$ is torsion-free rank-1 on 
$C(u)\x J^1(u)\big/J^1(u)$, and the proof is complete.

\sct3 Flag schemes

 \lem1 Let $X$ be a scheme, and $\cF$ a coherent $\cO_X$-module.  Assume
$\cF$ is invertible at each associated point of $X$, and is everywhere
locally generated by two sections.  Set $W:=\IP(\cF)$, and let
$w\:W\to X$ be the structure map.  Then Serre's graded
$\cO_X$-algebra homomorphism  $\alpha$ is an isomorphism:
        $$\alpha\:\cSym(\cF)\risom\bigoplus_{n\geq 0}w_*\cO_W(n).$$
 \pf The question is local on $X$, and $\cF$ is locally generated by two
sections.  So, setting $\cE:=\cO_X^{\oplus 2}$, we may assume there is a
short exact sequence of the form:
        $$0\to \cN\TO\ve \cE\TO\vf \cF\to 0.$$

This sequence induces a  short exact sequence of graded
$\cSym(\cE)$-modules:
        $$0\to\cN\ox\cSym(\cE)[-1]\to\cSym(\cE)\to
	\cSym(\cF)\to 0.\eqno\Cs1.1)$$
A priori, the sequence is only right exact.  However, since $\ve$ is
injective, every associated point of $\cN$ is an associated point $P$ of
$X$; so we need only check for left exactness at such a $P$.  By
hypothesis, $\cF$ is invertible at $P$; whence, $\cN$ is too.
Therefore, left exactness holds at $P$, so everywhere.

Set $V:=\IP(\cE)$, and let $v\:V\to X$ be the structure map.  Then $\vf$
induces a closed embedding $\iota\:W\into V$ such that $w=v\circ\iota$.
Moreover, applying the exact functor ``tilde''  to \Cs1.1), we obtain the
following short exact sequence on $V$:
        $$0\to v^*\cN\ox\cO_V(-1)\to\cO_V\to\cO_W\to 0,$$
 in which $\cO_V\to\cO_W$ is the comorphism of $\iota$.

For convenience, given a coherent  $\cO_X$-module  $\cG$ and integer
$i$, set  
	$$R(\cG,i):=R^1v_*(v^*\cG\ox\cO_V(i)),$$ 
 and let $b(\cG,i)$ denote the following natural map:
	$$b(\cG,i)\:\cG\ox\cSym_i(\cE)\To
	v_*(v^*\cG\ox\cO_V(i)).$$

For every $i\geq 0$, consider the natural commutative diagram:
        $$\CD 
        0@>>>\cN\ox\cSym_{i-1}(\cE) @>>> 
	\cSym_i(\cE) @>>> \cSym_i(\cF)@>>>0 \\
	@. @V{b(\cN,\,i-1)}VV @V{b(\cO_V,i)}VV @VVV \\
	0@>>>v_*(v^*\cN\ox\cO_X(i-1)) @>>> 
	v_*\cO_V(i) @>>> w_*\cO_W(i) @>>> R(\cN,i-1).
	\endCD$$
 Since $\cE$ is free, $b(\cO_V,i)$ is an isomorphism by Serre's
computation.  Therefore, to prove the lemma, it is
enough to prove that $R(\cN,\,i-1)=0$ and that $b(\cN,\,i-1)$ is an
isomorphism.

Fix $i\geq -1$.  Given a short exact sequence of coherent
$\cO_X$-modules
	$$0\to\cG'\to\cG\to\cG''\to 0,\eqno\Cs1.2)$$
 consider the following induced diagram:
        $$\CD
        0@>>> \cG'\ox\cSym_i(\cE) @>>> 
	\cG\ox\cSym_i(\cE) @>>>
	\cG''\ox\cSym_i(\cE)@>>>0  \\
	@. @Vb(\cG',i)VV @Vb(\cG,i)VV @Vb(\cG'',i)VV \\
	0@>>>v_*(v^*\cG'\ox\cO_V(i)) @>>> 
	v_*(v^*\cG\ox\cO_V(i)) @>>>
	v_*(v^*\cG''\ox\cO_V(i))@>>> R(\cG',i).
	\endCD$$
 Since $\cE$ is free, the upper sequence is exact.  Since $v$ is flat,
\Cs1.2) pulls back to a short exact sequence on $V$; whence, the lower
sequence is exact.

If $R(\cG',i)=0$ and if $b(\cG',i)$ and $b(\cG'',i)$ are isomorphisms,
then $b(\cG,i)$ is one too.  If, in addition, $R(\cG'',i)=0$, then also
$R(\cG,i)=0$.

On the other hand, if $b(\cG,i)$ and $b(\cG'',i)$ are isomorphisms, then
$b(\cG',i)$ is one too, and $R(\cG',i)\subseteq R(\cG,i)$.  If, in
addition, $R(\cG,i)=0$, then also $R(\cG',i)=0$.

Since $\cN\subset\cO_X^{\oplus 2}$, there is a 
short exact sequence
         $$0\to\cI\to\cN\to\cJ\to 0$$ 
 where $\cI\subseteq\cO_X$ and $\cJ\subseteq\cO_X$.  It is, therefore,
enough to prove that $R(\cG,i)=0$ and that $b(\cG,i)$ is an isomorphism
when $\cG\subseteq\cO_X$.

Finally, suppose  $\cG\subseteq\cO_X$, and let 
$Y\subseteq X$ be the subscheme defined by $\cG$.  
By Serre's computation, $b(\cO_X,i)$ and 
$b(\cO_Y,i)$ are isomorphisms.  Also, 
$R(\cO_X,i)=0$.  Hence, $R(\cG,i)=0$, and 
$b(\cG,i)$ is an isomorphism.  The proof is now complete.

\dsc2 Flag schemes.
 Let $f\:X\to T$ be a map of (locally Noetherian) schemes, and
$\cI$ a coherent sheaf on $X$.  As is conventional, for each $T$-scheme
$U$, set $X_U:=X\x U$, and let $\cI_U$ denote the pullback of $\cI$ to
$X_U$.  Fix $m\ge0$.  By a {\it $m$-flag of $\cI_U/X_U/U$}, let us mean
a chain of coherent sheaves,
	$$\cI_m\subset\cI_{m-1}\subset\cdots\subset\cI_1
	\subset\cI_0:=\cI_U,$$
 such that, for $1\le i\le m$, the $i$th quotient $\cI_{i-1}/\cI_i$ is
$U$-flat of relative length 1.  Denote the set of all these $m$-flags by
$\IF_m(U)$.

Since the quotients are flat, for each $U$-scheme $V$, the $m$-flag
pulls back to a $m$-flag of $\cI_V/X_V/V$.  So, as $U$ varies, the
$\IF_m(U)$ form a contravariant functor $\IF_m$.

Clearly, $\IF_0$ is representable by $T$.  Suppose $\IF_{m-1}$ is
representable by a $T$-scheme $F_{m-1}$, and consider the universal
$(m-1)$-flag:
	$$\cK_{m-1}\subset\cK_{m-2}\subset\cdots\subset
	\cK_1\subset\cK_0:=\cI_{F_{m-1}}.$$
Then, clearly, $\IF_m$ is representable by the Quot scheme
	$$F_m:=\Quot^1_{(\cK_{m-1}/X\x F_{m-1}/F_{m-1})};$$
furthermore, the universal  $m$-flag on $X\x F_m$  is the chain
	$$\cJ_m\subset\cJ_{m-1}\subset\cdots\subset
	\cJ_1\subset\cJ_0:=\cI_{F_m}$$
where $\cJ_i$ is the pullback of $\cK_i$ for $0\le i<m$ and where
$\cJ_m$ is the universal subsheaf of $\cJ_{m-1}$.  Call $F_m$ the {\it
$m$-flag scheme of $\cI/X/T$}.

According to \cite{K87, Prop.~(2.2), \p.~109}, we have $F_m =
\IP(\cK_{m-1})$; furthermore, if $(\phi_m,\tau_m)\:F_m\to X\x F_{m-1}$
denotes the structure map, then, on $X\x F_m$, we have
	$$\cJ_{m-1}/\cJ_m=(\phi_m,1)_*\cO_{F_m}(1).$$
 And there are natural maps,
	$$\CD 
	\psi^{(m,n)}\:
	F_m @>>> X\x F_{m-1} 
	@>>> X\x (X\x F_{m-2})
 	@>>>\cdots@>>> X^{\x n}\x F_{m-n}.
	\endCD$$  
 Set $\tau_{j,i}:=\tau_{i+1}\cdots\tau_j$; so $\tau_{j,i}\:F_j\to F_i$.
Set $\psi^{(m)}:=\psi^{(m,m)}$; so $\psi^{(m)}\:F_m\to X^{\x m}$.  Set
$\psi_i:= \phi_i\circ\tau_{m,i}$; so $\psi_i\:F_m\to X$.  Then
$\psi_m=\phi_m$ and $\psi^{(m)}= (\psi_m,\dots,\psi_1)$.  Let $q_i\:X\x
F_i\to F_i$ denote the projection.

\lem3 Under the conditions of \Cs2), assume that $f\:X\to T$ is a flat
projective family of integral curves that are locally planar, and assume
that $\cI$ is invertible.  Let $p_i$ and $p_{i,j}$ be the projections of
$X^{\x m}$ onto the indicated factors, and let $p$ be its structure map.
Let $\cX$ be the ideal of the diagonal $\Delta$  of $X\x X$.  Set
  $$\cH:=\psi^{(m)}_*\bigl(\cD_{q_m}(\cJ_m)^{-1}\bigr)
  \hbox{ and } \cF:=\cD_f(\cI)^{-1}$$
 where $\cD_{q_m}$ and $\cD_f$ mean determinant of cohomology.
Let $X_1,\dots,X_m\subset X$ be open subschemes such that 
$X_i\cap X_j\cap X_k$ is
$T$-smooth for all distinct $i,j,k$.  Set 
$H:=X_m\x\cdots\x X_1$.  Then
   $$\cH\mid H = \biggl(\biggl(\bigotimes_{i<j}p_{i,j}^*\cX\biggr)\ox
	\biggl(\bigotimes_i p_i^*\cI\biggr)\ox p^*\cF\biggr)\biggm| H.$$
 \pf  The additivity of the 
determinant of cohomology yields 
	$$\cD_{q_m}(\cJ_m)^{-1}
	=\biggl(\bigotimes_{i=1}^m
	\cD_{q_m}(\cJ_{i-1}/\cJ_i)\biggr)
	\ox\cD_{q_m}(\cI_{F_m})^{-1}.$$
Set $\cM_i:=\tau_{m,i}^*\cO_{F_i}(1)$.  Then
$\cJ_{i-1}/\cJ_i=(\phi_i,1)_*\cM_i$.  So $\cD_{q_m}(\cJ_{i-1}/\cJ_i)=\cM_i$.
Now, forming the determinant commutes with changing the base; so
	$$\cD_{q_m}(\cI_{F_m})^{-1}
	=\bigl(p\,\psi^{(m)}\bigr)^*\cD_f(\cI)^{-1}
	=\bigl(\psi^{(m)}\bigr)^*p^*\cF.$$
Since $\cF$ is invertible, we may apply the projection formula to
$p^*\cF$ and $\psi^{(m)}_*$.  We therefore have to prove that
      $$\psi^{(m)}_*\biggl(\bigotimes_{i=1}^m\cM_i\biggr)\biggm| H
	= \biggl(\biggl(\bigotimes_{i<j}p_{i,j}^*\cX\biggr) \ox
   \biggl(\bigotimes_i p_i^*\cI\biggr)\biggr)\biggm| H.\eqno\Cs3.1)$$

Since $F_1 = \IP(\cI)$ and $\cI$ is invertible, $F_1=X$ and $\cM_1=\cI$
and $\psi^{(m)}=1$.  So \Cs3.1) holds when $m=1$.  We proceed by
induction on $m$.

Suppose $m\ge2$.  Set $\cN_i:=\tau_{m-1,\,i}^*\cO_{F_i}(1)$ for
$i=1,\dots,m-1$.  Let $u_i$ and $u_{i,j}$ be the projections of
$X^{\x(m-1)}$ onto the indicated factors.  Set
	$$\cG:=\biggl(\bigotimes_{i<j}u_{i,j}^*\cX\biggr)\ox
	\biggl(\bigotimes_i u_i^*\cI\biggr)\and
	G:=X_{m-1}\x\cdots\x X_1.$$
Then the induction hypothesis yields
        $$\psi^{(m-1)}_*\biggl(\bigotimes_{i=1}^{m-1}\cN_i
	\biggr)\biggm|G=\cG\big|G.\eqno\Cs3.2)$$

Since $\cM_i= \tau_m^*\cN_i$ for  $i=1,\dots,m-1$, the 
projection formula yields
	$$(\phi_m,\tau_m)_*\biggl(\bigotimes_{i=1}^m\cM_i\biggr)
	=q_{m-1}^*\biggl(\bigotimes_{i=1}^{m-1}\cN_i\biggr)
	\ox (\phi_m,\tau_m)_*\cO_{F_m}(1).\eqno\Cs3.3)$$

For  $i=0,\dots,m-1$, set
        $$\cG_i:=\biggl(\bigotimes_{j=1}^ip_{1,m-j+1}^*\cX
	\biggr)\ox p_1^*\cI.$$ 
 For $i>0$, set $\wt\cG_{i-1}:=(u_{m-i},1)^*\cG_{i-1}$ and
$H_i:=(u_{m-i},1)^{-1}H$.  Now, $\cI$ is invertible.  Since
$u_{m-i}(H_i)\subseteq X_m\cap X_i$ and since $X_m\cap X_i\cap X_j$ is
$T$-smooth for each $j<i$, the pullback $\wt\cG_{i-1}$ is invertible.

Set $\wt G:=\bigl(\psi^{(m-1)}\bigr)^{-1}G\subset F_{m-1}$ and
$G_1:=X_m\x\wt G$.  For $i=1,\dots,m-1$, set
$\gamma_i:=1\x\psi_{m-1,i}$; so $\gamma_i\:X\x F_{m-1}\to X\x X$.  Let
us prove that
	$$\cK_i\big|G_1
   =\biggl(\biggl(\bigotimes_{j=1}^i\gamma_j^*\cX\biggr)
       \ox \cI_{F_{m-1}}\biggr)
	 \biggm| G_1\eqno\Cs3.4)$$
by induction on $i\ge0$.
 First off, $\cK_0=\cI_{F_{m-1}}$.  So suppose \Cs3.4) holds for $i-1$.
Let $\Gamma_i\subset X\x F_{m-1}$ be the image of $F_{m-1}$ under
$(\psi_{m-1,i},1)$.  Then $\Gamma_i=\gamma_i^{-1}\Delta$, where
$\Delta\subset X\x X$ is the diagonal.  Since $\Delta/X$ is flat,
forming its ideal $\cX$ commutes with changing the base.  Hence
$\gamma_i^*\cX$ is the ideal of $\Gamma_i$.

Set $\cO_{\Gamma_i}(1):=(\psi_{m-1,i},1)_*\cN_i$.  
Then $\cK_i$ is the kernel of the composition
	$$\cK_{i-1}\onto\cO_{\Gamma_i}\ox \cK_{i-1}
	\onto\cO_{\Gamma_i}(1).\eqno\Cs3.5)$$
 The kernel of the first surjection is the ideal-module product 
$(\gamma_i^*\cX)\cdot\cK_{i-1}$.  Hence
$(\gamma_i^*\cX)\cdot\cK_{i-1}\subset\cK_i$.  
Now, $\Gamma_i\cap G_1\subseteq (X_m\cap X_i)\x\wt G$.  
Also, $\cK_{i-1}\big|(X\x\wt G)$ is the structure sheaf, whence 
invertible, off the union of the $X_j\x\wt G$ for 
$j=1,\dots,i-1$. Since, by hypothesis, 
$X_m\cap X_i\cap X_j$ 
is $T$-smooth for  $j=1,\dots,i-1$, the restriction 
$\cK_{i-1}\big|(\Gamma_i\cap G_1)$ is invertible.  Thus the 
second surjection in \Cs3.5) is an isomorphism on $G_1$, and 
hence $(\gamma_i^*\cX\cdot\cK_{i-1})|G_1=\cK_i|G_1$.  
In addition, 
	$$(\gamma_i^*\cX\ox\cK_{i-1})|G_1\risom
	(\gamma_i^*\cX\cdot\cK_{i-1})|G_1,$$ 
 and therefore
	$$\cK_i\big|G_1=
	\bigl(\gamma_i^*\cX\ox\cK_{i-1}\bigr)\big| G_1.$$
 Combining the expression above with \Cs3.4) for $i-1$, 
we get \Cs3.4) for $i$.

Since $\psi^{(m-1)}=(\psi_{m-1,m-1},\dots,\psi_{m-1,1})$, 
Formula \Cs3.4) is equivalent to 
        $$\cK_i\big|G_1=((1\x\psi^{(m-1)})^*\cG_i)|G_1.
	\eqno\Cs3.6)$$

For $i=1,\dots,m-1$, set $\wt H_i:=(\psi_{m-1,i},1)^{-1}(G_1)$.  
As the second map in \Cs3.5) is an isomorphism on $G_1$, we have
	$(\psi_{m-1,i},1)^*\cK_{i-1}|\wt H_i= \cN_i|\wt H_i$.
 Now, note that  
	$$(1\x\psi^{(m-1)})\circ(\psi_{m-1,i},1)=
	 (u_{m-i},1)\circ\psi^{(m-1)};\eqno\Cs3.7)$$
 whence $(\psi^{(m-1)})^{-1}H_i=\wt H_i$.  Therefore, \Cs3.6) yields
	$$\cN_i|\wt H_i=(\psi^{(m-1)})^*
	\wt\cG_{i-1}|\wt H_i.\eqno\Cs3.8)$$

Consider a geometric point $x$ of a fiber of $H/G$, and view $x$ as well
as a point of $X_m$.  Suppose, at $x$, the ideal $p_{1,m-j+1}^*\cX$ is
not invertible.  Then $x\in X_m\cap X_j$.  And $p_{1,m-j+1}^*\cX$ is
generated by two elements owing to Nakayama's Lemma, since the fibers of
$X/T$ are locally planar.  Furthermore, $x\notin X_m\cap X_k$ for any
$k\neq j,m$ since $X_m\cap X_j\cap X_k$ is $T$-smooth by hypothesis;
whence $p_{1,m-k+1}^*\cX$ is invertible at $x$.  Now, $\cI$ is
invertible.  Therefore, it follows from \Cs3.6) that $\cK_i|G_1$ is
everywhere locally generated by two sections.

Since $F_m=\IP(\cK_{m-1})$, Lemma \Cs1) yields
	$$(\phi_m,\tau_m)_*\cO_{F_m}(1)\big|G_1
	 = \cK_{m-1}\big|G_1.$$
Therefore, it follows from \Cs3.3) that
        $$(\phi_m,\tau_m)_*\biggl(\bigotimes_{i=1}^m\cM_i\biggr)
	\biggm|G_1=
	\biggl(q_{m-1}^*\biggl(\bigotimes_{i=1}^{m-1}\cN_i
	\biggr)\ox\cK_{m-1}\biggr)\biggm|G_1.$$

Let $r\:X^{\x m}\to X^{\x m-1}$ be the projection onto the product of
the last factors.  To complete the proof, it is now enough to prove that
        $$(1\x\psi^{(m-1)})_*
	\biggl(q_{m-1}^*\biggl(\bigotimes_{i=1}^{m-1}\cN_i
	\biggr)\ox\cK_k\biggr)\biggm|H=
	(r^*\cG\ox\cG_k)\bigm|H\eqno\Cs3.9)$$
for $k=0,\dots,m-1$.  Again, we proceed by induction.

 First off, $\cK_0=\cI_{F_{m-1}}$, and $\cI$ is invertible.  So the
projection formula yields
        $$(1\x\psi^{(m-1)})_*
	\biggl(q_{m-1}^*\biggl(\bigotimes_{i=1}^{m-1}\cN_i
	\biggr)\ox\cK_0\biggr)\biggm|H=\biggl(p_1^*\cI\ox
	r^*\psi^{(m-1)}_*
	\biggl(\bigotimes_{i=1}^{m-1}\cN_i\biggr)
	\biggr)\biggm|H.$$
But   $\cG_0:=p_1^*\cI$.  Thus \Cs3.2) yields \Cs3.9) when $k=0$.

Suppose \Cs3.9) holds for $k-1$ with $k<m$.  Now, owing to the
discussion surrounding \Cs3.5), the exact sequence
	$$0\to\cG_k\big|H\to\cG_{k-1}\big|H\to
	(\cG_{k-1}\ox p_{1,m-k+1}^*\cO_{\Delta})\big|H\to 0$$
 pulls back under $1\x\psi^{(m-1)}$ to the exact sequence
         $$0\to\cK_k\big|G_1\to\cK_{k-1}\big|G_1\to
	 \cO_{\Gamma_k}(1)\big|G_1\to 0.$$
Hence \Cs3.9) for $k$ follows from \Cs3.9) for $k-1$ provided
         $$(1\x\psi^{(m-1)})_*
	\biggl(q_{m-1}^*\biggl(\bigotimes_{i=1}^{m-1}\cN_i
	\biggr)\ox\cO_{\Gamma_k}(1)\biggr)
	\biggm|H=
	(r^*\cG\ox\cG_{k-1}\ox p_{1,m-k+1}^*\cO_{\Delta})\bigm| H.$$

Here, the right-hand side is equal to
$(u_{m-k},1)_*(\cG\otimes\wt\cG_{k-1})\big|H$.  On the other hand, since 
$\cO_{\Gamma_k}(1):=(\psi_{m-1,k},1)_*\cN_k$, 
the left-hand side is equal to
         $$(1\x\psi^{(m-1)})_*(\psi_{m-1,k},1)_*
	 \biggl(\biggl(\bigotimes_{i=1}^{m-1}\cN_i
	\biggr)\ox\cN_k\biggr)\biggm|H,$$
 or because of \Cs3.7), to
	$$(u_{m-k},1)_*\psi^{(m-1)}_*
	\biggl(\biggl(\bigotimes_{i=1}^{m-1}\cN_i
	\biggr)\ox\cN_k\biggr)\biggm|H.$$

Finally, observe that
	$$\psi^{(m-1)}_*\biggl(\biggl(\bigotimes_{i=1}^{m-1}\cN_i
	\biggr)\ox\cN_k\biggr)\biggm|H_k=
	(\cG\otimes\wt\cG_{k-1})\big|H_k;$$
Indeed, this formula results from \Cs3.8) and from the projection
formula, since, as noted above, $\wt\cG_{k-1}\big|H_k$ is invertible.  The
proof is now complete.

\sct4 Isomorphism

\thm1 Let $C/S$ be a flat projective family of integral curves with, at
worst, ordinary nodes and cusps.  Then the autoduality
isomorphism $\smash{\Pic^0_{\CJ/S}}\risom J$ extends uniquely to an
isomorphism of compactifications $\Pic^-_{\CJ/S}\risom\CJ$.
 \pf Set $U:=\smash{\Pic^0_{\CJ/S}}$ and $\?U:=\smash{\Pic^-_{\CJ/S}}$.
By Theorem (2.6), the autoduality isomorphism extends uniquely to a map,
say $\eta\:\?U\to\CJ$.  By descent theory, it suffices to prove $\eta$
is an isomorphism after a faithfully flat base change.  So we may assume
the smooth locus of $C/S$ admits a section $\sigma$.  Set
$\cL:=\cO_C(\sigma(S))$.  Then the autoduality isomorphism is simply
$A_\cL^*$.  Also, $C\x\CJ/\CJ$ carries a universal sheaf $\cI$, which is
of degree 0 and rigidified along $\sigma$.  Finally, we may assume as
well that $S$ is the spectrum of a Henselian local ring with
algebraically closed residue field.

Let $\beta\:J\to U$ be the inverse of $\smash{A_\cL^*}$, and set
$U^=:=\smash{\Pic^=_{\CJ/S}}$.  It suffices to extend $\beta$ to a map
$\?\beta\:\CJ\to U^=$.  Indeed, $\CJ$ is the schematic closure of $J$;
so $\?\beta$ factors through $\?U\subseteq U^=$.  Also,
$\eta\?\beta\big|J=1$, so $\eta\?\beta=1$.  And $\?\beta\eta\big|U=1$;
so $\?\beta\eta=1$.  Hence $\eta$ is an isomorphism.  Thus it suffices
to construct $\?\beta\:\CJ\to U^=$.

First, recall from \cite{EGK02, Prop.~2.2, \p.595}, or rather from
\cite{EGK02, Rmk.~2.4, \p.597}, the definition of $\beta\:J\to U$.  Let
$q_1$ and $q_2$ be the projections of $C\x\CJ$ onto the indicated
factors.  Form the product $\cK:=\cI\ox q_1^*\cL^{\ox (g-1)}$ on $C\x\CJ$
and its determinant of cohomology $\cP:=\cD_{q_2}(\cK)$ on $\CJ$.  Let
$r_1\:\CJ\x J\to\CJ$ be the projection, and $\mu\:\CJ\x J\to\CJ$ the
multiplication map.  Set
	$$\cT:=r_1^*\cP\ox\mu^*\cP^{-1},$$
Then $\cT$ is invertible and defines $\beta$.

(Note that there are two canceling sign errors in \cite{EGK02, Rmk.~2.4,
\p.597}: first, the theta divisor is the zero scheme of the canonical
regular section of $\cP^{-1}$, not $\cP$; second, there $\cT$ is the
inverse of what it should be.  With these corrections, the discussion in
\cite{EGK02, Rmk.~2.4} goes through.)

Let $C_0\subseteq C$ be the smooth locus of $C/S$.  Set $Z:=C-C_0$ and
equip $Z$ with its induced reduced subscheme structure.  Let $s\in S$ be
the closed point, and $P_1,\dots,P_m\in C(s)$ the singularities.  For
$i=1,\dots,m$, there exists, by \cite{EGA, IV$_4$ 18.5.11, \p.130}, a
decomposition $Z=Z_i\cup\wt Z_i$ in which $Z_i$ and $\wt Z_i$ are
disjoint closed subschemes such that one of them, say $Z_i$, meets
$C(s)$ only in $P_i$.

Set $C_i:=C-\wt Z_i$ for $i=1,\dots,m$.  Then $C_i$ is open, and
$C_i\supset C_0$.  Also $P_i\in C_i$.  So the closed fiber $C(s)$ is
covered by the $C_i$; whence, $C=C_1\cup\cdots\cup C_m$.  If $i\neq j$,
then $C_i\cap C_j$ does not contain any point of $Z$ in the closed fiber
$C(s)$; whence, $C_i\cap C_j\subseteq C_0$.  Therefore, $C_i\cap
C_j=C_0$.

View the bigraded Abel map as a rational map $A\:C\x\CJ\to\CJ^{-1}$.
Set $\nu:=\mu_\cL\circ A$.  Viewing $C^{\x (i+1)}$ as $C^i\x C$, form
$(1\times\nu)\: C^{\x (i+1)}\x\CJ\to C^i\x\CJ$ and
       $$C^{\x m}\x\CJ\to C^{\x m-1}\x\CJ\to\cdots\to 
	C\x\CJ\to\CJ.$$
By \cite{EGK00, Cor. 2.6, \p.5969}, the rational map $A$ is smooth where
$A$ is defined; whence, the composition is smooth where it is defined.
Set
	$$H_m:=C_m\x\cdots\x C_2\x C_1\subset C^{\x m}
	\and W:=H_m\x J.$$
Then the composition is defined on $W$ because $C_i\cap C_j=C_0$ if
$i\neq j$.  Thus there is a well-defined smooth map $\zeta\:W\to\CJ$.
Again since $C_i\cap C_j=C_0$ if $i\neq j$, it follows that
$\zeta^{-1}J=C_0^{\x m}\x J$.

Not only is $\zeta\:W\to\CJ$ smooth, but also $\zeta$ is surjective.
Indeed, since $P_1,\dots,P_m$ are ordinary nodes or cusps, every
torsion-free rank-1 sheaf on $C(s)$ is of the form $\cJ\ox\cL$ where
$\cL$ is invertible and $\cJ$ is the ideal of a reduced subscheme of
$\bigcup P_i$;  hence, $\zeta(W)\supset \CJ(s)$.
But $\zeta$ is smooth, so open.  Therefore, $\zeta(W)=\CJ$.

Below, we construct a torsion-free rank-1 sheaf $\cQ$ on $\CJ\x W\big/W$
such that $\cQ$ coincides with $(1\x\zeta)^*\mu^*\cP^{-1}$ on $\CJ\x
C_0^{\x m}\x J$.  Using $\cQ$, we can complete the proof as follows.
First, let $\? r_1\:\CJ\x W\to\CJ$ denote the projection, and set
$\cT_\zeta:=\? r_1^*\cP\otimes\cQ$.  Then $\cT_\zeta$ is also a torsion-free
rank-1 sheaf on  $\CJ\x W\big/W$.  So $\cT_\zeta$ induces a
map $\smash{\?\beta'}\:W\to U^=$.  But $\cT_\zeta$ and $(1\x\zeta)^*\cT$
coincide on $\CJ\x C_0^{\x m}\x J$.  Hence $\smash{\?\beta'}$ and
$\beta\circ\zeta$ coincide on $C_0^{\x m}\x J$.

Set $V:=W\x_\CJ W$, and let $\zeta_1$ and $\zeta_2$ be the projections.
Set $V_0:=\zeta^{-1}_1\zeta^{-1}(J)$.  Then
$V_0=\zeta^{-1}_2\zeta^{-1}(J)$.  So $\?\beta'\circ\zeta_1\big|V_0 =
\?\beta'\circ\zeta_2\big|V_0$.  Now, $\zeta$ is flat, and
$J\supset\Ass(\CJ)$; so $V_0\supset\Ass(V)$.  Hence
$\?\beta'\circ\zeta_1=\?\beta'\circ\zeta_2$.  Therefore, descent theory
yields a unique map $\?\beta\:\CJ\to U^=$ such that
$\?\beta'=\?\beta\circ\zeta$.  Since
	$$\?\beta\circ\zeta\bigm|C_0^{\x m}\x J
	=\?\beta'\bigm|C_0^{\x m}\x J=
	\beta\circ\zeta\bigm|C_0^{\x m}\x J,$$ 
 and $\zeta$ is faithfully flat, $\?\beta|J=\beta$.  In other words,
$\?\beta$ is the desired extension of $\beta$. 

It remains to construct $\cQ$.  To lighten the notation, given
$S$-schemes $U$ and $V$ and coherent sheaves $\cG_1$ and $\cG_2$ on
$C\x U$ and $C\x V$, denote by $\cG_1\bt\cG_2$ the tensor product on
$C\x U\x V$ of the pulled-back sheaves.  Let $\cX$ denote the ideal of
the diagonal of $C\x C$, and for each $n>0$ and $i=1,\ldots,n$,  set
	$$\cE^{(n)}_i:=p_{1,n-i+2,\ldots,n+1}^*
	\bigl(\cX\bt\cdots\bt\cX\bigr)$$
 where $p_{1,n-i+2,\ldots,n+1}\:C^{n+1}\to C^{i+1}$ is the projection
onto the indicated factors.

Given any $n\ge0$, let $F_n$ be the $n$-flag scheme of $\cI/(C\x\CJ)/\CJ$,
and let
	$$\cI^{(n)}_n\subset\cI^{(n)}_{n-1}
	\subset\cdots\subset\cI^{(n)}_1\subset
	\cI^{(n)}_0:=\cI_{F_n}$$
 be the universal flag.  Form the natural map 
$\psi^{(n)}\:F_n\to C^{\x n}\x\CJ$, and set
	$$F'_n:=(\psi^{(n)})^{-1}(C_0^{\x n}\x\CJ).$$
Also, let $\iota\:C^{\x n}\x\CJ\to\CJ\x C^{\x n}$ be the switch map, and
set $\wh\psi^{(n)}:=\iota\circ\psi^{(n)}$.

Note that $\psi^{(n)}\big|F'_n$ is an isomorphism, whose inverse is 
defined by the $n$-flag
   $$(\cE^{(n)}_n|(C\x C_0^{\x n}))\bt\cI\subset\cdots\subset
   (\cE^{(n)}_1|(C\x C_0^{\x n}))\bt\cI\subset\cI_{C_0^{\x n}\x\CJ};$$
this flag is well-defined since  $\cX$ is invertible on $C_0\x C_0$. 

Set $\cN:=\cI|(C\x J)$.  Then $\cI^{(n)}_n\bt\cN\bt\cL^{\ox n}$ is a
torsion-free rank-1 sheaf on $C\x F_n\x J\big/F_n\x J$; so it defines a
map $\gamma_n\:F_n\x J\to\CJ$.

Take $n:=m$.  Set $\wt F_m:=(\psi^{(m)})^{-1}(H_m\x\CJ)$ and 
$\wt F'_m:=(\psi^{(m)})^{-1}(W)$.  Note that 
$\psi^{(m)}\big|\wt F'_m$ is an isomorphism, whose inverse is 
defined by the $m$-flag
   $$(\cE_m^{(m)}|(C\x H_m))\bt\cN\subset\cdots\subset
   (\cE^{(m)}_1|(C\x H_m))\bt\cN\subset\cI_W;$$
this flag is well-defined since $\cN$ is invertible and $C_i\cap
C_j=C_0$ if $i\neq j$.

 Next, let us see that the following square is commutative:
	$$\CD
              F'_m\x J           @>\gamma_m>> \CJ \\
          @V\wh\psi^{(m)}\x 1VV                   @A\mu AA    \\
        \CJ\x C_0^{\x m}\x J  @>1\x\zeta >>  \CJ\x J.
	\endCD$$
 Indeed, here $\wh\psi^{(m)}\x 1$ is an isomorphism; in addition,
$\gamma_m\circ(\wh\psi^{(m)}\x 1)^{-1}$ and $\mu\circ(1\x\zeta)$ coincide,
because they are defined by the same sheaf, namely,
	$$\cI\bt(\cE^{(m)}_m|(C\x C_0^{\x m}))
      \bt\cN\bt\cL^{\ox m} \hbox{ on } C\x\CJ\x C_0^{\x m}\x J.$$

Since the square is commutative and $\wh\psi^{(m)}\x 1$ is an isomorphism,
 $(1\x\zeta)^*\mu^*\cP^{-1}$ and $(\wh\psi^{(m)}\x 1)_*\gamma_m^*\cP^{-1}$
coincide on $\CJ\x C_0^{\x m}\x J\subset \CJ\x W$.  So set
	$$\cQ:=(\wh\psi^{(m)}\x 1)_*\gamma_m^*\cP^{-1}\bigm|\CJ\x W.
	\eqno\Cs1.1)$$
 It remains to show that $\cQ$ is torsion-free rank-1 on $\CJ\x W\big/W$.
Since $\zeta$ is faithfully flat, it suffices to form the map $\zeta\x 1
\: W\x W\to \CJ\x W$ and show that $(\zeta\x 1)^*\cQ$ is torsion-free
rank-1 on  $W\x W\big/ W$. 

Given any $p,q\ge0$, take $n:=p+q$, and let us construct the
following square:
	$$\CD
	F_n     @>\theta>>   F_p		\\
   @V\tau_{n,q}VV         @VV\tau_{p,0}V	\\
        F_q	@>\omega>>	\CJ.
	\endCD\eqno\Cs1.2)$$
Here  $\tau_{n,q}$ and $\tau_{p,0}$  are the maps defined in (3.2).  Let
$\theta$ be  defined by the $p$-flag
  $$\cI^{(n)}_{n}\bt\cL^{\ox q}\subset\cdots
  \subset\cI^{(n)}_{q}\bt\cL^{\ox q}$$
on $C\x F_{n}\x \CJ\big/F_{n}\x \CJ$; here $\cI^{(n)}_{q}\bt\cL^{\ox q}$
is torsion-free, rank-1, and of degree 0 since $\cI^{(n)}_0$ is so,
since each quotient $\cI^{(n)}_{i-1}/\cI^{(n)}_i$ is flat of relative
length 1, and since $\cL$ is invertible of degree 1.  Let $\omega$ be
defined by the torsion-free, rank-1, and degree-0 sheaf
$\cI^{(q)}_q\bt\cL^{\ox q}$.  Clearly, the square is Cartesian.

Note that the following square is commutative:
	$$\CD
	F_n     	@>\theta>>   	F_p		\\
   @V\psi^{(n,p)}VV         	   @VV\psi^{(p)}V	\\
        C^{\x p}\x F_q	@>1\x \omega>>	C^{\x p}\x\CJ.
	\endCD\eqno\Cs1.3)$$
It follows formally that this square is Cartesian since \Cs1.2) is so.

Take $p=q:=m$.  Recall that $\psi^{(m)}$ restricts to the isomorphism
$\wt F'_m\risom W$.  Set 
$\wt F_{2m}:=(\psi^{(2m)})^{-1}(H_m\x W)$.  
Then Square \Cs1.3) yields this one:
	$$\CD
	\wt F_{2m}     	@>\theta>>   	\wt F_m		\\
   @V\psi^{(2m)}VV         	   @VV\psi^{(m)}V	\\
        H_m\x W		@>1\x\zeta >>	H_m\x\CJ.
	\endCD\eqno\Cs1.4)$$

Interchange the two factors in the lower left and right corners of 
Square \Cs1.4), and multiply on the right by $J$.  
The result is the following Cartesian square:
	$$\CD
 \wt F_{2m}\x J		@>\theta\x1>>	\wt F_m\x J		\\
     @V\xi VV				@VV\wh\psi^{(m)}\x1V	\\
        W\x W		@>\zeta\x1>>	\CJ\x W,
	\endCD\eqno\Cs1.5)$$
which introduces $\xi$. In fact, $F_{2m}\x J$ can be viewed as the 
$2m$-flag scheme of $(\cI\bt\cN)/C\x\CJ\x J/\CJ\x J$, 
and $\xi$ is, up to interchanging factors, the restriction of the natural 
map $F_{2m}\x J\to C^{\x 2m}\x\CJ\x J$ over $H_m\x H_m\x J\x J$.

Since the lower map $\zeta\x1$ in Square \Cs1.5) is flat, Equation
\Cs1.1) yields
	$$(\zeta\x1)^*\cQ = \xi_*(\theta\x1)^*\gamma_m^*\cP^{-1}.$$
So it remains to show $\xi_*(\theta\x1)^*\gamma_m^*\cP^{-1}$ is
torsion-free rank-1 on $W\x W\big/ W$.

Recall  $\cP:=\cD_{q_2}(\cK)$ where  $\cK:=\cI\bt\cL^{\ox (g-1)}$.
Now, on  $C\x F_{2m}\x J$,  we have
	$$(1\x\theta\x1)^*(1\x\gamma_m)^*\cK
	=\cI^{(2m)}_{2m}\bt\cN\bt\cL^{\ox 2m}.$$
 Let $q_{2,3}\:C\x F_{2m}\x J\to F_{2m}\x J$ be the projection.  Since
forming the determinant of cohomology commutes with changing the base,
it remains to see that
 	$$\xi_*D_{q_{2,3}}\bigl(\cI^{(2m)}_{2m}\bt\cN\bt\cL^{\ox 2m}\bigr)
	\bigm|\bigl(\wt F_{2m}\x J\bigr)$$
is torsion-free rank-1 on $W\x W\big/ W$.  But it is clearly so in
view of Lemma (3.3) applied with 
$C\x J\x J/J\x J$ for $X/T$, with $\cN\bt\cN\bt\cL^{\ox 2m}$ for $\cI$, and 
with $X_i:=C_i\x J\x J$ for $i=1,\dots,m$ and
$X_i:=C_{i-m}\x J\x J$ for $i=m+1,\dots,2m$.

 \references
 \serial{duke}{Duke Math. J.}
 \serial{zeit}{Math. Zeitschrift}
 \serial{jag}{J. Algebraic Geom.}
 \serial{adv}{Adv. Math.}
 \serial{ajm}{Amer. J. Math.}
 \serial{bams}{Bull. Amer. Math. Soc.}
 \serial{ca}{Comm. Alg.}
 \serial{jlms}{J. London. Math. Soc. (2)}
 \serial{bpm}{Birkh\"auser Prog. Math.}
 \def\splnm #1 #2 {Springer LNM {\bf#1}\ (#2)}
 \def\GrothFest#1 #2 {\unskip, \rm
 in ``The Grothendieck Festschrift. Vol.~#1,'' P.~Cartier, 
 L.~Illusie, N.M.~Katz, G.~Laumon, Y.~Manin, K.A.~Ribet (eds.) 
 \bpm 86 1990 #2 }
 \def\osloii{\unskip, \rm in ``Real and complex singularities. 
 Proceedings, Oslo 1976,'' P. Holm (ed.), Sijthoff \& Noordhoff, 
 1977, pp.~}
 \def\sitgesii{\unskip, \rm in ``Enumerative Geometry. Proceedings,
 Sitges, 1987,'' S. Xamb\'o Descamps (ed.), \splnm 1436 1990 }

AIK76
 {A. Altman, A. Iarrobino, and S. Kleiman},
 Irreducibility of the compactified Jacobian
 \osloii 1--12.

AK76
 A. Altman and S. Kleiman,
 Compactifying the Jacobian,
 \bams 82(6) 1976 947--49

AK79
 A. Altman and S. Kleiman,
 Compactifying the Picard scheme II,
 \ajm 101 1979 10--41

AK80
 A. Altman and S. Kleiman,
 Compactifying the Picard scheme,
 \adv 35 1980 50--112

EGK00
 E. Esteves{,} M. Gagn\'e{,} and S. Kleiman,
 Abel maps and presentation schemes,
 \ca 28(12) 2000 5961--5992

EGK02
 E. Esteves{,} M. Gagn\'e{,} and S. Kleiman,
 Autoduality of the compactified Jacobian,
(preliminary version {\tt math.AG/9911071}),
\jlms 65 2002 591--610

FGA
 A. Grothendieck,
 Technique de descente et th\'eor\`emes d'existence en g\'eom\'etrie
 alg\'e\-bri\-que V, S\'eminaire Bourbaki, 232, Feb. 1962, and VI,
 S\'eminaire Bourbaki, 236, May 1962.

EGA
 A. Grothendieck and J. Dieudonn\'e,
 El\'ements de G\'eom\'etrie Alg\'ebrique,
 Publ. Math. IHES, IV$_4$, Vol. 32, 1967.

K87
 S. Kleiman,
 Multiple-point formulas II: the Hilbert scheme
 \sitgesii , 101--138.

L02
 G. Laumon,
 Fibres de Springer et jacobiennes compactifi\'ees,
 {\tt math.AG/0204109.}

S04
 J. Sawon,
 Derived equivalence of holomorphic symplectic manifolds,
 {\tt math.AG/0404365.}

\endreferences
\bye